%% file: main.tex
\documentclass[11 pt]{amsart}
\usepackage{amssymb}
\usepackage{amsmath}
\usepackage{amsthm}
\usepackage{amsfonts}
\usepackage{standalone}

\usepackage[margin=1 in]{geometry}
\usepackage{multicol}
\usepackage{tikz}
\usetikzlibrary{patterns}
\usepackage{lmodern,tikz}
\usetikzlibrary{positioning}
\usepackage{pgfplots,caption}
\pgfplotsset{compat=1.13}
\usepackage[utf8]{inputenc}
\usepackage{graphicx}
\usepackage{verbatim}
\usepackage{enumerate}
\usepackage{mathtools}
\usepackage{xcolor}
\usetikzlibrary{calc}
\usepackage{float}

\usepackage{mathpazo}
\usepackage{euler}
\usepackage{microtype}

\DeclarePairedDelimiter\abs{\lvert}{\rvert}%
\DeclarePairedDelimiter\norm{\lVert}{\rVert}%

\makeatletter
\let\oldabs\abs
\def\abs{\@ifstar{\oldabs}{\oldabs*}}
\let\oldnorm\norm
\def\norm{\@ifstar{\oldnorm}{\oldnorm*}}
\makeatother

\usepackage{amsmath,stackengine,scalerel}
\stackMath
\def\dhat#1{\ThisStyle{\setbox0=\hbox{$\SavedStyle#1$}%
  \stackengine{0pt}{\SavedStyle#1}{\SavedStyle\hspace{.4\ht2}%
  \hat{\vphantom{#1}}\kern\dimexpr3.2\LMpt+2.0pt\relax\hat{\vphantom{#1}}}{O}{c}{F}{T}{L}}%
}

\tikzset{
    triple/.style args={[#1] in [#2] in [#3]}{
        #1,preaction={preaction={draw,#3},draw,#2}
    }
}

\newcommand{\LE}{%
 \hbox{%
  \vbox{\hrule width 0.35em height 0.04 em}%
  \vbox{\offinterlineskip%
   \hbox{\kern -0.02em\vrule height 0.65em width 0.04em \hspace{0.1em}}%
  }%
 }%
}

\DeclareMathAlphabet{\mathcal}{OMS}{cmsy}{m}{n}

\newcommand{\M}{\mathcal{M}}

\newcommand{\newword}[1]{\textbf{\emph{#1}}}

\usepackage{xpatch}
\makeatletter
\xpatchcmd{\@thm}{\thm@headpunct{.}}{\thm@headpunct{}}{}{}
\makeatother
\newtheorem{thm}{Theorem}[section]
\newtheorem{theorem}[thm]{Theorem}

\newtheorem{proposition}[thm]{Proposition}
\newtheorem{lemma}[thm]{Lemma}

\newtheorem{conjecture}[thm]{Conjecture}

\newtheorem{corollary}[thm]{Corollary}

\newtheorem{definition}[thm]{Definition}

\newtheorem{remark}[thm]{Remark}

\newtheoremstyle{case}{}{}{}{}{}{:}{ }{}
\theoremstyle{case}

\tikzset{
    single/.style args={}{
        draw,line width=.5mm,black
    }
}
\tikzset{
    double/.style args={}{
        preaction={draw,line width=1.25mm,black},draw,line width=.25mm,white
    }
}
\tikzset{
    triple/.style args={}{
        preaction={preaction={draw,line width=2mm,black},draw,line width=1mm,white},draw,line 
        width=.5mm,black
    }
}  
\tikzset{
    quadruple/.style args={}{
        preaction={preaction={preaction={draw,line width=2.75mm,black},draw,line width=1.75mm,white},draw,line width=1.25mm,black},draw,line width=.25mm,white
    }
}

\tikzset{
    unmarked/.style args={}{
        draw, fill=black, circle, inner sep=0pt, minimum size=.25cm,
    }
}

\tikzset{
    marked/.style args={}{
        draw, fill=white, circle, inner sep=0pt, minimum size=0.25cm, line width=0.06cm
    }
}

\long\def\basic{
        \fill [white] (-0.3, -0.3) rectangle (2.3, 3.3);

        \node [unmarked={}] (1) at (0.5,0.5) {};
        \node [unmarked={}] (2) at (0.5,1.5) {};
        \node [unmarked={}] (3) at (0.5,2.5) {};
        \node [unmarked={}] (6) at (1.5,1.5) {};
        \node [unmarked={}] (7) at (1.5,2.5) {};
        
        \node [unmarked={}] (9) at (2. ,2.5) {};
        \node [unmarked={}] (8) at (2. ,1.5) {};
        \node [unmarked={}] (5) at (1.5,1. ) {};
        \node [unmarked={}] (4) at (1. ,0.5) {};
        \node [unmarked={}] (0) at (0.5,0. ) {};
        
        \draw [line width = .3mm, black](0,0) -- (1,0);
        \draw [line width = .3mm, black](0,1) -- (2,1);
        \draw [line width = .3mm, black](0,2) -- (2,2);
        \draw [line width = .3mm, black](0,3) -- (2,3);

        \draw [line width = .3mm, black](0,0) -- (0,3);
        \draw [line width = .3mm, black](1,0) -- (1,3);
        \draw [line width = .3mm, black](2,1) -- (2,3);
}

\def\xSep{2}
\def\ySep{4.5}

\long\def\bork#1#2#3{
    \node [marked={}] at (#1) {};
    \node [marked={}] at (#2) {};
    \node [marked={}] at (#3) {};
}

\title{The $h$-vector of a Positroid is a pure O-sequence}
\author{Amy He, Pierce Lai, SuHo Oh}

\address{Bellaire Senior High School}
\email{he2023amy@gmail.com}

\address{Massachusetts Institute of Technology, Department of Electrical Engineering and Computer Science}
\email{pwlai@mit.edu}

\address{Texas State University, Department of Mathematics}
\email{s\_o79@txstate.edu}

\date{\today}

\begin{document}


\begin{abstract}
A well-known conjecture of Stanley is that the $h$-vector of any matroid is a pure ${\mathcal O}$-sequence. There have been numerous papers with partial progress on this conjecture, but it is still wide open. Positroids are special class of linear matroids that play a crucial role in the field of total positivity. In this short note, we prove that Stanley's conjecture holds for positroids.
\end{abstract}

\maketitle

\section{Introduction}
A matroid is a structure which generalizes the notion of independence that rises in linear algebra, graph theory and other areas. From a matroid, we can construct a simplicial complex (\newword{matroid complex}) and study its \newword{h-vector} which provides the topological information of the complex. The following conjecture by Stanley is still wide open:

\begin{conjecture}[\cite{Stanley1977}]
\label{con:stanley}
The h-vector of  a matroid is a pure O-sequence.
\end{conjecture}


There are some specific classes of matroids for which the above conjecture has been established to be true. In particular, these include cographic matroids by Merino in \cite{Merino2001}, lattice-path matroids by Schweig in \cite{Schweig2010}, cotransversal matroids by Oh in \cite{Oh2013}, paving matroids by Merino, Noble, Ramirez-Ibanez, and Villarroel-Flores  \cite{MNRV2012}, internally perfect matroids by Dall in \cite{Dal}, rank $3$ matroids by H\'a, Stokes, and Zanello in \cite{HaStokesZanello2013}, rank $3$ and corank $2$ matroids by DeLoera, Kemper, and Klee in \cite{DelKemKle}, rank $4$ matroids by Klee and Samper in \cite{KleeSamper2015}, rank $d$ matroids with $h_d \leq 5$  by Constantinescu, Kahle, and Varbaro in \cite{ConKahVar} and some small special classes of graphic matroids by Kook \cite{Kook2011}, Preston et al \cite{biconed} and David et al \cite{triconed}.

In this paper the focus will solely be on the class of positroids, which is a special class of realizable matroids, for which the conjecture is still wide open. Positroids come from the study of the nonnegative part of the Grassmannian \cite{Postnikov} and have seen increased applications in physics, with use in the study of scattering amplitudes \cite{arkani} and the study of shallow water waves \cite{kodama}. Positroids contain the class of lattice-path matroids, but is incomparable with all other classes mentioned above. The closest class is perhaps that of transversal matroids since positroids and transversal matroids are both contained in the class called gammoids. Regardless, the class of positroids and the class of transversal matroids (and also the class of cotransversal matroids since dual of a positroid is again a positroid) are independent \cite{marcott2019basis}.

\section{Preliminary}

Throughout the paper, we will use the following notation. When $A$ is a set and $a$ is a single element, instead of writing $A \setminus \{a\}$, we will  write $A \setminus a$. Similarly, instead of writing $A \cup \{a\}$, we will  use $A \cup a$.

\subsection{Positroids} 

\begin{figure}
\begin{multicols}{2}

 \begin{center}
 
  \vspace{3.7mm}
    \input{tikzpics/lediagram}
    \captionsetup{width=1.0\linewidth}
  \captionof{figure}{A $\LE$-diagram.}
  \label{fig:lediagram}

 \end{center}

\vfill\null
\columnbreak

 \begin{center}
 
  \vspace{3.7mm}
    \input{tikzpics/lenetwork}
    \captionsetup{width=1.0\linewidth}
  \captionof{figure}{The corresponding $\LE$-network.}
  \label{fig:lenetwork}

 \end{center}
\end{multicols}
\end{figure}

We will start with the definition of $\LE$-diagrams and $\LE$-graphs and use that to define positroids. A $\LE$-diagram will be a Young diagram with filled dots in some boxes. Consider the southeastern boundary path of the diagram. We call this the \newword{boundary path} of the diagram. Starting from the northeast corner and heading towards the southwest corner, label each edge of the boundary path with integers $1,\ldots,n$. We call these the \newword{boundary labels}. We index each column and row of a diagram with the label of the unique boundary label it contains. A box at $(i,j)$ stands for the box of $L$ at row indexed with $i$ and column indexed with $j$. For example, a box of $L$ in Figure~\ref{fig:lediagram} at $(2,5)$ is the box at the second row from top (indexed with $2$) and first column from left (indexed with $5$). We write $L_{i,j} = 1$ if there the box at $(i,j)$ contains a dot, and $L_{i,j} = 0$ otherwise.

\begin{definition}[\cite{Postnikov}, Definition 6.1] For a partition $\lambda$, let us define a $\LE$-diagram $L$ of shape $\lambda$ as a filling of boxes of the Young diagram of shape $\lambda$ such that, for any three boxes at $(i,j),(i',j),(i,j')$, where $i'<i$ and $j'>j$, if $L_{i',j} = L_{i,j'} = 1$ then $L_{i,j} = 1$. This property is called the $\LE$-property.

\end{definition}

From each $\LE$-diagram $L$, we can get a planar network called the $\LE$-graph in Definition 6.3 of \cite{Postnikov}, via the following way: place a dot in the middle of each step in the boundary lattice path of the diagram and mark these dots with the corresponding boundary labels (so $1,\ldots,n$ from top to bottom). We will call these vertices the \newword{boundary vertices}. For each dot inside the $\LE$-diagram, draw a horizontal line to its right, and vertical line to its bottom until it reaches another dot. Next orient all vertical edges downward and horizontal edges to the right.  The vertical labels will be the source vertices and the horizontal labels will be the sink vertices in the $\LE$-graph. Since each $\LE$-diagram encodes a unique $\LE$-graph, we will abuse notation and go between a $\LE$-diagram and its corresponding $\LE$-graph freely in this paper. An example of a $\LE$-graph obtained from the $\LE$-diagram in Figure~\ref{fig:lediagram} is drawn in Figure~\ref{fig:lenetwork}.

The set $[n]$ (where $n$ is the length of the boundary path) is called the \newword{ground set} (being consistent with the notation for the positroid that will be indexed by the $\LE$-diagram). Let $B_0$ denote the set of source vertices in the $\LE$-graph of $L$. For example, in the $\LE$-diagram $L$ drawn in Figure~\ref{fig:lediagram}, we have $B_0 = 124$ as we can see from the corresponding $\LE$-graph in Figure~\ref{fig:lenetwork}. A \newword{path} in a $\LE$-graph is a directed path that starts at some boundary vertex and ends at some boundary vertex. A \newword{family of non-touching paths} is a family of paths where no two paths intersect within the graph. We say that a family of non-touching paths consisting of paths $p_1,\ldots,p_k$ \newword{represents} a set $B_0 \setminus \{p_1^s,\ldots,p_k^s\} \cup \{p_1^e,\ldots,p_k^e\}$ where $p_i^s$ stands for the starting point and $p_i^e$ stands for the endpoint of the path $p_i$. For example, the bottom picture in Figure~\ref{fig:allthebases} shows an empty family, which represents the set $124$. The leftmost picture at the top row in Figure~\ref{fig:allthebases} shows a family consisting of two paths, one starting at $1$ and ending at $3$, the other one starting at $4$ and ending at $5$: this represents the set $124 \setminus 14 \cup 35 = 235$. Notice that any set represented by a family of non-touching paths in a $\LE$-graph will have the same cardinality as $B_0$.


From Theorem~6.5 of \cite{Postnikov}, there is a bijection between positroids and $\LE$-diagrams. Moreover, the same theorem gives the following result, which we will use as the definition of a positroid in this paper:

\begin{definition}
\label{def:pos}
Fix a positroid $\M_L$ that corresponds to a $\LE$-diagram $L$. Then the set of bases of $\M_L$ is exactly the collection of sets represented by a family of non-touching paths in the $\LE$-graph of $L$.
\end{definition}

From the above definition, the set $B_0$ is naturally the lexicographic minimal basis of $\M_L$. We call this the \newword{canonical basis} of the positroid. For example, take a look at the $\LE$-diagram in Figure~\ref{fig:lediagram}. All the bases we can get are:
$$124,125,134,135,145,234,235,245,345,$$ 
as one can check from various families of disjoint paths drawn on the $\LE$-graph as in Figure~\ref{fig:allthebases}. 

An element $e$ of the ground set is called a \newword{loop} if it is not contained in any basis and is called a \newword{coloop} if it is contained in every basis.

\subsection{Activity and the h-vector} 
From any ordering on the ground set $E$, we get an induced lexicographic ordering on the bases. An element $i$ of a basis $B$ is \newword{(internally) active} if $B \setminus i \cup j$ is not a base for any $j < i$. Otherwise, it is called \newword{(internally) passive}. There is also the notion of externally active and externally passive elements, but in this paper we will only stick to interally active and passive elements, and hence will skip the usage of the word internally throughout. Given a basis $B$ of a matroid, we call its \newword{passivity} as the number of (internally) passive elements it contains.

We will use this activity property of h-vectors as the definition of h-vectors of matroids in this paper.

\begin{theorem}
Let $(h_0,\dots,h_r)$ be the $h$-vector of a matroid $\M$. For $0 \leq i \leq r$, the entry $h_i$ is the number of bases of $\M$ with $i$ passive elements.
\end{theorem}

Let us again look at the running example of Figure~\ref{fig:allthebases}. Take a look at $135$. Here $1$ is active since it cannot be replaced with a smaller element to get a basis of $\M$. On the other hand, $3$ and $5$ are passive since we can replace either of them with $2$ to get a different basis of $\M$. The $h$-vector of the corresponding positroid is $(1,2,3,3)$.

\begin{remark}
\label{rem:noloopcoloop}
The class of positroids is closed under deletion and contraction \cite{ARW}. Deleting a loop (simply remove that element from the ground set) or contracting by a coloop (remove that element from the ground set and all the bases) does not change the $h$-vector of the matroid. For sake of proving Stanley's conjecture, we may assume that our matroid does not have any loops or coloops. Hence, from now on, throughout the entire paper, we are only going to deal with positroids with no loops nor coloops. In terms of the $\LE$-diagram, it means there is no empty row nor column. 
\end{remark}

\subsection{Order ideal and pure O-sequences}
An \newword{order ideal} is a finite collection $X$ of monomials such that, whenever $M \in X$ and $N$ divides $M$, we have $N \in X$ as well. We define a maximal monomial as the monomial such that all of other monomials in the order ideal are divisors of the maximal monomial. Given an order ideal its \newword{O-sequence} is a sequence that can be obtained from counting monomials of each degree. We say that the order ideal is \newword{pure} if all the maximal monomials are of the same degree, and also call the O-sequence pure if it comes from a pure order ideal.

For example, take a look at Figure~\ref{fig:allthebases} and focus only on the monomials. If we take any monomial among this collection, say $x_3x_5^2$, we can see that all divisors of it is again present in the collection. Hence this is an order ideal. Moreover, since the maximal monomials under division (which are usually not guaranteed to all have the same degree) here are $x_3^2x_5,x_3x_5^2,x_5^3$ which have the same degree, and hence this is a pure order ideal. Notice that the degree sequence of this pure order ideal is $(1,2,3,3)$. Hence we can see that the $h$-vector of the positroid coming from Figure~\ref{fig:lediagram}, which was $(1,2,3,3)$, is indeed a pure $O$-sequence. We have just verified that Stanley's conjecture is true for this case. What we want to do in this paper is to similarly construct a pure order ideal for each and every positroid, hence proving the conjecture for this class.

\section{The main result}

Our goal in this section is to come up with an algorithm that maps a family of non-crossing paths to a monomial, using the elements of $[n] \setminus B_0$ as the variables. Then we will show that the image under this map is a pure order ideal to answer Stanley's conjecture for positroids.

Given a path from a source $s$ to sink $k$, we write it as $(s \rightarrow k)$. Since each basis $B \in \M$ is coming from a family of non-crossing paths from set of sources to set of sinks, we write $B$ as a collection of $(s \rightarrow k)$'s and call this the \newword{non-crossing matching representation} (or simply a matching representation) of $B \in \M$. Note that all sources $s$ are elements of $B_0 \setminus B$, and all sinks $k$ are elements of $B \setminus B_0$. 

\begin{figure}[h]
 \begin{center}
 
  \vspace{3.7mm}
    \input{tikzpics/basis/allthebases}
    \captionsetup{width=1.0\linewidth}
  \captionof{figure}{Hasse diagram of the $\LE$-Diagrams which form the bases of the matroid. Each $\LE$-Diagram is labelled with its corresponding monomial.}
  \label{fig:allthebases}

\end{center}
\end{figure}
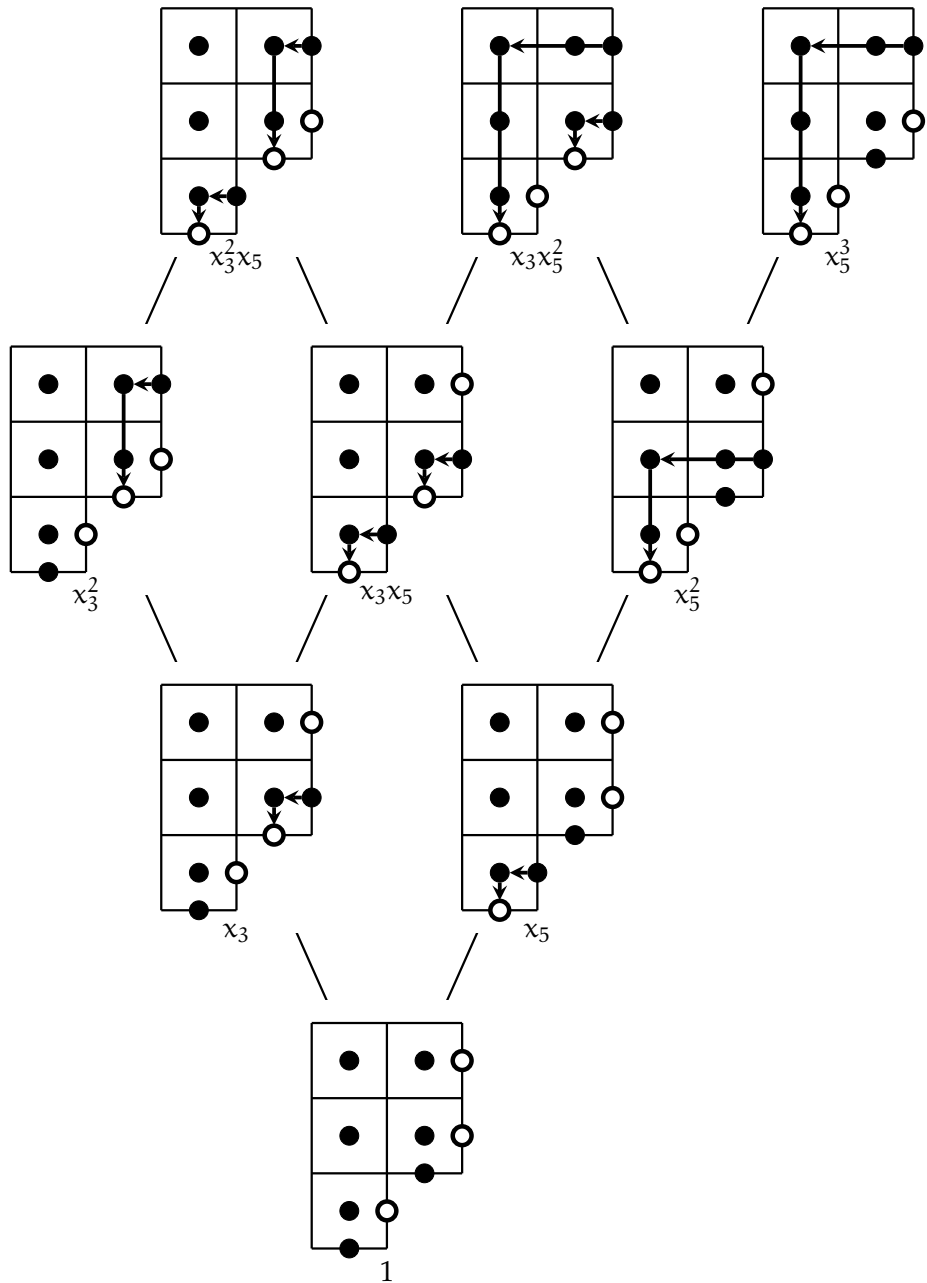

We will slightly abuse notation and go between the path $(s \rightarrow k)$ and the corresponding interval $(s,k)$ freely. So for example, the path $(2 \rightarrow 5)$ \newword{contains} $3$ since $3 \in (2,5)$. Given a matching representation $R$, we put a partial ordering $\prec_R$ on elements of $B_0 \setminus B$ (and also separately on elements of $B \setminus B_0)$ by setting $s'  \prec  s$ ($k' \prec k$ respectively) if $(s \rightarrow k),(s' \rightarrow k') \in R$ and $s < s' < k' < k$: i.e. if the path coming out of $s'$ is nested inside (as intervals on $[n]$) the path coming out of $s$. For example, in the middle picture of the top row of Figure~\ref{fig:allthebases}, we have $2 \prec_R 1$ (also $3 \prec_R 5)$ since comparing the paths $2 \rightarrow 3$ and $1 \rightarrow 5$, the interval $(2,3)$ is nested inside the interval $(1,5)$. We put an ordering on the paths so it is compatible with the ordering we just discussed on the starting point: for example we have $(2,3) \prec_R (1,5)$. 

\begin{lemma}
\label{lem:uniquerep}
Fix a positroid $\M_L$. For each basis $B \in \M_L$, there exists a unique matching representation of $B$.
\end{lemma}

\begin{proof}
From Definition~\ref{def:pos}, there exists at least one matching representation of each $B \in M_L$. For sake of contradication, assume that for some $B \in M_L$ we have two different matching representations $R_1$ and $R_2$. Consider the partial order $\prec_{R_1}$ on $B_0 \setminus B$.  Under this ordering, pick any minimal element $s$ of $B_0 \setminus B$ such that $(s \rightarrow k) \in R_1$ and $(s \rightarrow k') \in R_2$ with $k \not = k'$. If $k < k'$, then pick $s' \in B_0 \setminus B$ so that $(s' \rightarrow k) \in R_2$. Since $R_2$ is non-crossing we must have $s < s'$. But then from the way $s$ was chosen, we need to have $(s' \rightarrow k) \in R_1$ as well and we get a contradiction. If $k > k'$, then pick $s' \in B_0 \setminus B$ so that $(s' \rightarrow k') \in R_1$. Since $R_1$ is non-crossing we must have $s < s' < k' < k$. But then from the way $s$ was chosen, we must have $(s' \rightarrow k') \in R_2$ and we again get a contradiction.
\end{proof}

It is worth noting that although we are calling $s \rightarrow k$ a path, there might be multiple different paths in the $\LE$-graph that are encoded by the same $s \rightarrow k$. For example in the $\LE$-graph of Figure~\ref{fig:lediagram}, consider the path that starts from $1$, goes left, left, down, down, down to $5$. This path is encoded as $1 \rightarrow 5$ in the matching representation. If we look at the path that starts from $1$, goes left, down, left, down, down to $5$, this is a different path from above but is still encoded as $1 \rightarrow 5$. So even though Lemma~\ref{lem:uniquerep} tells us that each $B$ has a unique matching representation, there may be a number of different families of non-crossing paths that correspond to that same matching representation.

Thanks to Lemma~\ref{lem:uniquerep}, when considering the partial order $\prec_R$ coming from a matching representation of some basis $B \in \M_L$, we may instead write it as $\prec_B$. Using this, we slightly abuse notation and say $(s \rightarrow k) \in B$ if $(s \rightarrow k)$ is a path inside the unique matching representation of $B$. Given some fixed basis $B \in \M_L$, we say that an element $x \in [n]$ is \newword{associated under $B$} (or simply associated, when it is clear we are dealing with the basis $B$) to a path $(s \rightarrow k) \in B$ if $(s \rightarrow k)$ is the unique minimal path containing $x$ under the ordering $\prec_B$. Next we show that each passive element $b$ in a basis of $B$ is closely related to the matching representation of $B$.

\begin{lemma}
\label{lem:containment}
Fix a positroid $\M_L$ and a basis $B \in \M_L$. For each passive $b \in B\cap B_0$, there is a path $(s \rightarrow k) \in B$ containing $b$ (i.e. $s < b < k$).
\end{lemma}

\begin{proof}
Pick $b \in B \cap B_0$ that is passive. Assume for sake of contradiction there is no path $(s \rightarrow k) \in B$ such that $s < b < k$. From $b$ being passive, there is some $B' \in \M_L$ such that $B' = B \setminus b \cup b'$ with $b' < b$ and $b' \not \in B$. We have $(b \rightarrow k) \in B'$ for some $k \in B' \setminus B_0$. Then there is some $c$ such that $(c \rightarrow k) \in B$. From the assumption we have $c > b$. 

Under the ordering $\prec_{B'}$, pick any minimal element $x$ such that $(x \rightarrow y) \in B'$ and $(x' \rightarrow y) \in B$ with $x' > x \geq b$. Such $x$ exists thanks to the argument in the previous paragraph. Since $(b,n] \cap B = (b,n] \cap B'$, we have $x' \not \in B'$. So there is some $y'$ such that $(x' \rightarrow y') \in B'$, with $x <x'<y'<y$. Then, this implies $y'\in B$, so there exists some $x'' > x'$ such that $(x'' \rightarrow y')\in B$. Since $x' \prec_{B'} x$, we get a contradiction with $x$ being chosen as a minimal element within $\prec_B'$.


\end{proof}

We are now ready to present the map $\phi$ which will be mapping each basis of a positroid into a monomial supported on variables indexed by $[n] \setminus B_0$. Given some basis $B \in M_L$, for each element $b \in B \setminus B_0$, we simply map it to the variable $x_b$. For a passive element $b \in B \cap B_0$, we look at the unique minimum (under $\prec_B$) path $s \rightarrow k$ that contains $b$ (recall that we called this as $b$ being associated to $s \rightarrow k$ under $B$), then map $b$ to $x_k$. Now take the product of all such variables to obtain a monomial, which is going to be supported on variables indexed by $[n] \setminus B_0$ as we wanted. 

For example, take a look at the leftmost picture in the top row of Figure~\ref{fig:allthebases}. This family of paths encodes the basis $235$. The passive elements of this basis is actually all of $2,3,5$. Since $3$ and $5$ are sinks of paths, $3$ gives $x_3$ and $5$ gives $x_5$. Now the minimum (under $\prec_{235}$) path containing $2$ is $(1 \rightarrow 5)$ so $2$ gives $x_3$. The resulting monomial is $x_3^2 x_5$.

\begin{proposition}
\label{prop:phi}
Let $\M_L$ be a positroid. Given a basis $B$, we have $\phi(B) = \prod_{k \in B \setminus B_0} x_k^{|T_k^B|+1}$ where $T_k^B$ consists of elements $x$ associated to $(s \rightarrow k) \in B$ such that replacing $s \rightarrow k$ with $x \rightarrow k$ gives some other basis in $\M_L$.
\end{proposition}
\begin{proof}
From the way $\phi$ is defined, it is enough to show that $T_k^B$ consists exactly of the passive elements of $B \cap B_0$ associated to  $(s \rightarrow k) \in B$. It is clear that if we can replace $s \rightarrow k$ with $x \rightarrow k$ for $x \in (s,k)$ then $x$ is passive. Hence we only need to show the other direction: if $x$ is passive then we can replace $s \rightarrow k \in B$ with $x \rightarrow k$ to get another basis.

Since $x$ is passive, there is some $x' < x$ such that $B' = B \setminus x \cup x' \in \M_L$. Consider the matching representation of $B$. Since we have a path $(s \rightarrow k) \in B$, from the non-crossing property we have equal number of source points and sink points within $(s,k)$. Moreover, since $x$ is associated to $s \rightarrow k$, we have equal number of source points and sink points within $(s,x)$ and also within $(x,k)$. This property for $(x,k)$ carries over to the matching representation of $B'$ as well. From this property, we have $(x \rightarrow k) \in B'$.
\end{proof}

For example, again take a look at the leftmost picture in the top row of Figure~\ref{fig:allthebases}. The set $T_5^{235}$ is empty. The set $T_3^{235}$ consists of $2$, since we have $2 \in (1,3)$ and we can replace $1 \rightarrow 3$ with $2 \rightarrow 3$ to get the set $135$, which is again a basis of the positroid.

The Proposition~\ref{prop:phi} together with Lemma~\ref{lem:containment} tells us that $\phi$ maps the passivity of the basis to the degree of the monomial:

\begin{corollary}
Let $\M_L$ be positroid. Then the passivity of $B$ is equal to the degree of $\phi(B)$.
\end{corollary}

For example, take a look at the leftmost picture of second row from the top within Figure~\ref{fig:allthebases}. We have the path $1 \rightarrow 3$ so the corresponding basis is $234$. Here $2$ and $3$ are passive, but $4$ is active. So $234$ has passive degree to be $2$. The corresponding monomial is $x_3^2$, which again has degree of $2$.

\begin{lemma}
The map $\phi$ is a one-to-one map.
\end{lemma}

\begin{proof}
Assume for sake of contradiction that we can find two different bases $B, B'$ of the positroid $\M_L$ such that $\phi(B) = \phi(B')$. In $B$, under the ordering $\prec_B$, pick the minimal path $(s \rightarrow k)$ that isn't contained in $B'$. Since $\phi(B) = \phi(B')$, there has to be some path $(s' \rightarrow k) \in B'$ with $s \not = s'$. If $s < s'$, from the way $s \rightarrow k$ was chosen, we can replace $s \rightarrow k$ with $s' \rightarrow k$ within $B$ to get a basis so we have $T_k^{B'} \subsetneq T_k^{B}$. If $s > s'$, we cannot have a path in $B'$ that starts within $(s',s)$ and ends within $(s,k)$ due to the way $s \rightarrow k$ was chosen. This means we can replace $s' \rightarrow k$ with $s \rightarrow k$ within $B'$ to get a basis so we have $T_k^{B'} \supsetneq T_k^{B}$.
\end{proof}

\begin{proposition}
\label{prop:poi}
The map $\phi$ gives a pure order ideal.
\end{proposition}
\begin{proof}
We first start with the order ideal property. We need to show that given any basis $B$, for each variable $x_k$ appearing in $\phi(B)$, we can find a basis $B'$ such that $\phi(B') = \phi(B)/x_k$. To use induction, we introduce a new statistic on the poset with elements $B \setminus B_0$ and ordering $\prec_B$: we set $\sigma(k)$ to stand for the length of the longest possible chain having $k$ as the lower end within the poset. 

We are going to use induction on $\sigma(k)$. We start with the base case, when $\sigma(k) = 1$, that is when the path $(s \rightarrow k) \in B$ is not contained in any other path. If $|T_k^B| = 1$, we simply remove the path to get the $B'$ we want. Otherwise, we replace $s$ with the smallest element of $|T_k^B|$ to get our desired $B'$.

Now assume for sake of induction we have proven the result for $\sigma(k) < t$ and are now considering the case $\sigma(k) = t$. Same as above, if $|T_k^B| = 1$, we simply remove the path to get the $B'$ we want. Otherwise, we replace $s$ with the smallest element of $|T_k^B|$ to get some other basis $B''$. Then $\phi(B'')$ is obtained from $\phi(B)$ by dividing by $x_k$, but the powers of other variables $x_{k'}$ might have gone up for $k \prec_B k'$. Regardless, in such cases we have $\sigma(k') < t$, so we can find a new basis by shedding them off thanks to the induction assumption. This gives us the $B'$ we were looking for.

Now for the pure part, fix a basis $B$. We know that every element of $B \setminus B_0$ is automatically passive since we can simply undo the corresponding path to get a smaller basis. So potential active elements only occur within $B \cap B_0$. We want to show that if there is some active element, we can always find some basis $B'$ such that $\phi(B')$ is divisible by $\phi(B)$. Let $a$ be the largest element of $B$ that is active. We do a case by case analysis based on whether $a+1$ is in $B_0$ or not.

When $a+1 \not \in B_0$: from Remark~\ref{rem:noloopcoloop} we have $L_{a,a+1} = 1$. If $a+1 \in B$, then we have $s \rightarrow a+1$ in $B$ which is replaceable with $a \rightarrow a+1$ to get $B'$. This means $a$ is passive so we can ignore this case and say $a+1 \not \in B$. In this case, we may add the path $a \rightarrow a+1$ to get $B'$ we want ($\phi(B') = \phi(B) x_{a+1})$. 

When $a+1 \in B_0$: first assume for sake of contradiction that $a+1 \in B$. Then it is passive thanks to the way $a$ was chosen. This means we have $s \rightarrow k \in B$ such that $a+1 \in T_k^B$. From the $\LE$-property, the path $s \rightarrow k$ is replaceable with $a \rightarrow k$ as well, and we get a contradiction from Proposition~\ref{prop:phi}. 
    
Hence we have that $a+1 \in B_0 \setminus B$. This means we have some $(a+1 \rightarrow k) \in B$. Let $k'$ be the minimal element such that $a \rightarrow k'$ is possible inside the $\LE$-graph. If $k' \leq k$, from the $\LE$-property we can switch $a+1 \rightarrow k$ to $a \rightarrow k$ to get the $B'$ we desire ($\phi(B') = \phi(B) x_{a+1}$). In the case $k' > k$, from the way $k'$ was chosen and the $\LE$-property, we have $L_{a,k'} = 1$ and $L_{i,j} = 0$ for $i \leq a$ and $j \leq k'$ except when $(i,j) = (a,k')$. If the dot at $(a,k')$ is already occupied in any family of non-crossing paths representing $B$, we can replace the starting point of this path with $a$. This would contradict $a$ being active, so we may assume that the dot at $(a,k')$ is never occupied. If there is some path $s \rightarrow k' \in B$ with $s' > a$, then replace $s \rightarrow k'$ with $a \rightarrow k'$ to get the desired $B'$ ($\phi(B') = \phi(B)x_{k'}^t$ with $t \geq 1$). Otherwise we can simply add $a \rightarrow k'$ to get the $B'$ we desire ($\phi(B') = \phi(B) x_{k'}$).

\end{proof}

Thanks to the above proposition, we have established Stanley's conjecture for the class of positroids:

\begin{corollary}
The $h$-vector of any positroid is a pure $O$-sequence.
\end{corollary}

An example is shown in Figure~\ref{fig:allthebases}. Using the $\LE$-diagram from Figure~\ref{fig:lediagram}, all the bases and the corresponding family of non-crossing paths are drawn in Figure~\ref{fig:allthebases}. The monomials from the mapping $\phi$ are written below each picture. Proposition~\ref{prop:poi} guarantees that the collection of monomials we have obtained is indeed a pure order ideal.

\section*{Acknowledgements}
The research started under the $2021$-High school math camp hosted at Texas State University. The authors would like to thank the camp organizers for providing support and a great working environment.

\bibliography{bibl}
\bibliographystyle{siam}

\end{document}

%% file: tikzpics/lediagram.tex
    \begin{tikzpicture}[scale=1.5]

        \node [unmarked={}] (1) at (0.5,0.5) {};
       \node [unmarked={}] (2) at (0.5,1.5) {};
        \node [unmarked={}] (3) at (0.5,2.5) {};
        \node [unmarked={}] (6) at (1.5,1.5) {};
        \node [unmarked={}] (7) at (1.5,2.5) {};
        
            
        
        \draw [line width = .3mm, black](0,0) -- (1,0);
        \draw [line width = .3mm, black](0,1) -- (2,1);
        \draw [line width = .3mm, black](0,2) -- (2,2);
        \draw [line width = .3mm, black](0,3) -- (2,3);

        \draw [line width = .3mm, black](0,0) -- (0,3);
        \draw [line width = .3mm, black](1,0) -- (1,3);
        \draw [line width = .3mm, black](2,1) -- (2,3);
        
    \end{tikzpicture}

%% file: tikzpics/lenetwork.tex
    \begin{tikzpicture}[scale=1.5]

        \node [unmarked={}] (1) at (0.5,0.5) {};
        \node [unmarked={}] (2) at (0.5,1.5) {};
        \node [unmarked={}] (3) at (0.5,2.5) {};
        \node [unmarked={}] (6) at (1.5,1.5) {};
        \node [unmarked={}] (7) at (1.5,2.5) {};
        
        \node [unmarked={},label=right:{\large1}] (9) at (2. ,2.5) {};
        \node [unmarked={},label=right:{\large2}] (8) at (2. ,1.5) {};
        \node [unmarked={},label=right:{\large4}] (4) at (1. ,0.5) {};
        \node [unmarked={},label=below:{\large3}] (5) at (1.5,1. ) {};
        \node [unmarked={},label=below:{\large5}] (0) at (0.5,0. ) {};
            
        \draw [stealth-, line width = .5mm, black](0) -- (3);  
        \draw [stealth-, line width = .5mm, black](3) -- (9);  
        \draw [stealth-, line width = .5mm, black](2) -- (8);  
        \draw [stealth-, line width = .5mm, black](5) -- (7);  
        \draw [stealth-, line width = .5mm, black](1) -- (4);  
        
        \draw [line width = .3mm, black](0,0) -- (1,0);
        \draw [line width = .3mm, black](0,1) -- (2,1);
        \draw [line width = .3mm, black](0,2) -- (2,2);
        \draw [line width = .3mm, black](0,3) -- (2,3);

        \draw [line width = .3mm, black](0,0) -- (0,3);
        \draw [line width = .3mm, black](1,0) -- (1,3);
        \draw [line width = .3mm, black](2,1) -- (2,3);
        
    \end{tikzpicture}

%% file: tikzpics/basis/allthebases.tex
    \begin{tikzpicture}[scale=1]
    
        \begin{scope}[shift={(1,1.5)}]
            \draw [line width = .3mm] (0,0) -- (-\xSep,\ySep) -- (-2*\xSep,2*\ySep) -- (-1*\xSep,3*\ySep) -- (0,2*\ySep) -- (\xSep,\ySep);        
    
            \draw [line width = .3mm] (\xSep,3*\ySep) -- (2*\xSep,2*\ySep);
            \draw [line width = .3mm] (0,0) -- (3*\xSep,3*\ySep);
            \draw [line width = .3mm] (-\xSep, \ySep) -- (\xSep,3*\ySep);
        \end{scope}
    
        \begin{scope}[shift={(0, 0)}]
            \basic
            \bork{9}{8}{4}
            
            \node[] at (1,-0.3) {1};
        \end{scope}
        
        \begin{scope}[shift={(-\xSep,\ySep)}]
            \basic
            \bork{9}{5}{4}
            
            \draw [-stealth, line width = .5mm] (8) -- (6);
            \draw [-stealth, line width = .5mm] (6) -- (5);
            \node[] at (1,-0.3) {$x_3$};
        \end{scope}
        
        \begin{scope}[shift={(\xSep,\ySep)}]
            \basic
            \bork{9}{8}{0}
            
            \draw [-stealth, line width = .5mm] (4) -- (1);
            \draw [-stealth, line width = .5mm] (1) -- (0);
            \node[] at (1,-0.3) {$x_5$};
        \end{scope}
        
        \begin{scope}[shift={(-2*\xSep,2*\ySep)}]
            \basic
            \bork{8}{5}{4}
            
            \draw [-stealth, line width = .5mm] (9) -- (7);
            \draw [-stealth, line width = .5mm] (7) -- (5);
            \node[] at (1,-0.3) {$x_3^2$};
        \end{scope}
        
        \begin{scope}[shift={(0,2*\ySep)}]
            \basic
            \bork{9}{5}{0}
            
            \draw [-stealth, line width = .5mm] (8) -- (6);
            \draw [-stealth, line width = .5mm] (6) -- (5);
            \draw [-stealth, line width = .5mm] (4) -- (1);
            \draw [-stealth, line width = .5mm] (1) -- (0);
            \node[] at (1,-0.3) {$x_3x_5$};
        \end{scope}
        
        \begin{scope}[shift={(2*\xSep,2*\ySep)}]
            \basic
            \bork{9}{4}{0}
            
            \draw [-stealth, line width = .5mm] (8) -- (2);
            \draw [-stealth, line width = .5mm] (2) -- (0);
            \node[] at (1,-0.3) {$x_5^2$};
        \end{scope}
        
        \begin{scope}[shift={(-1*\xSep,3*\ySep)}]
            \basic
            \bork{8}{5}{0}
            
            \draw [-stealth, line width = .5mm] (9) -- (7);
            \draw [-stealth, line width = .5mm] (7) -- (5);
            \draw [-stealth, line width = .5mm] (4) -- (1);
            \draw [-stealth, line width = .5mm] (1) -- (0);
            \node[] at (1,-0.3) {$x_3^2x_5$};
        \end{scope}
        
        \begin{scope}[shift={(1*\xSep,3*\ySep)}]
            \basic
            \bork{5}{4}{0}
            
            \draw [-stealth, line width = .5mm] (9) -- (3);
            \draw [-stealth, line width = .5mm] (3) -- (0);
            \draw [-stealth, line width = .5mm] (8) -- (6);
            \draw [-stealth, line width = .5mm] (6) -- (5);
            \node[] at (1,-0.3) {$x_3x_5^2$};
        \end{scope}
        
        \begin{scope}[shift={(3*\xSep,3*\ySep)}]
            \basic
            \bork{8}{4}{0}
            
            \draw [-stealth, line width = .5mm] (9) -- (3);
            \draw [-stealth, line width = .5mm] (3) -- (0);
            \node[] at (1,-0.3) {$x_5^3$};
        \end{scope}

    \end{tikzpicture}